\documentclass[12pt,a4paper]{amsart}

\usepackage{amsmath,amsfonts,amssymb,mathrsfs}
\usepackage{graphicx}

\newcommand*{\leftstar}[1]{\vphantom{#1}_\star #1}
\newcommand*{\act}{\!\cdot\!}
\newcommand*{\init}{{\mathsf i}}
\newcommand{\mult}[2]{\multicolumn{2}{c|}{#1\act #2}}

\DeclareMathOperator{\card}{{\mathrm card}}

\newtheorem{definition}{Definition}

\begin{document}
\title{Moore machines duality}
\author[Jacques Peyri\`ere]{Jacques Peyri\`ere}

\address{Renmin University of China, School of Mathematics, 100872 Beijing, P.R.\ China, and
  Institut de Mathématiques d’Orsay, CNRS, Université Paris-Saclay,
  91405 Orsay, France.}

\email{jacques.peyriere@umiversite-paris-saclay.fr}

\date{June 22, 2020}
\maketitle

\begin{abstract} We present a simple algorithm to find the Moore
  machine with the minimun number of states equivalent to a given one.
\end{abstract}

\section{Introduction}

It is known that there is a unique, up to isomorphism, deterministic
finite automaton with a minimum number of states which recognizes the
same language as a given one. As explained in~\cite{hopcroft} (pages
29--30), the states of this minimal automaton can be taken as the
classes modulo an equivalence relation on the input words.

Several algorithms to minimize a deterministic finite automaton has
been proposed~\cite{bonchi,brz,moore}. According to~\cite{berstel}
these algorithms fall within three categories. Some of them apply as
well to non deterministic finite automata and to Moore machines.

In Section~\ref{algo} we describe an algorithm to minimize a Moore
machine, and establish its consistency. This algorithm lies in the
same family as Brzozowski's one~\cite{brz}. But its description as
well as the proof of its correctness are very simple and do not
require much apparatus. Furthermore, this gives another proof, simple
and elementary, of the existence and uniqueness of a minimum automaton
(see Section~3). In Section~4 we extend these results to morphisms of
free monoids.

\section{Moore machine and duality}

\begin{definition} A Moore machine is a 6-tuple
  $\mathscr M =(Q,\Sigma,\Delta,\delta,{\lambda},\init)$ where
  \begin{itemize}
  \item [--] $Q$ is a finite set, the set of states,
  \item [--] $\Sigma$ is a finite set, the input alphabet,
  \item [--] $\Delta$ is a finite set, the output alphabet,
  \item [--] $\delta$ is a mapping from $Q\times \Sigma$ to~$Q$,
    called the transition function,
  \item [--] ${\lambda}\in \Delta^Q$ is a mapping from~$Q$ to~$\Delta$,
  \item [--] $\init\in Q$ is the initial state.
  \end{itemize}
\end{definition}

\subsection{Actions of words on $Q$}

The set of words on the alphabet~$\Sigma$, including the empty
word~$\epsilon$, is denoted by~$\Sigma^*$. The concatenation of the
words~$u$ and~$v$ is simply written $uv$.\medskip

The transition function~$\delta$ can be extended in two ways as a
mapping from $Q\times \Sigma^*$ to~$Q$.

\begin{itemize}
\item[]\textsl{\large Right action.} We define, by recursion on the
  length of~$w$ a mapping, written $(a,w)\mapsto a\act w$, from
  $Q\times \Sigma^*$ to~$Q$:
  \begin{itemize}
  \item for $a\in Q$, $a\act \epsilon = a$,
  \item for $a\in Q$ and $j\in \Sigma$, $a\act j= \delta(a,j)$,
  \item $a\act wj = (a\act w)\act j$ for $a\in Q,\ w\in \Sigma^*$, and
    $j\in \Sigma$.
  \end{itemize}
\item[] \textsl{\large Left action.} In the same way we define a mapping,
  written $(w,a)\mapsto w\act a$, from $\Sigma^*\times Q$ to~$Q$:
  \begin{itemize}
  \item for $a\in Q$, $\epsilon\act a = a$,
  \item for $a\in Q$ and $j\in \Sigma$, $j\act a= \delta(a,j)$,
  \item $jw\act a = j\act(w\act a)$ for $a\in Q,\ w\in \Sigma^*$, and
    $j\in \Sigma$.
  \end{itemize}
\end{itemize}

\subsection{Right and left machines}

We can feed ${\mathscr M}$ with a word in~$\Sigma^*$ on the right or
on the left and get an output:
$${\mathscr M}\act w = {\lambda}({\mathsf i}\act w) \text{\quad and \quad}
w\act{\mathscr M} = {\lambda}(w\act {\mathsf i}).$$
If a state does not belong to the set
$\{\init\act w\ :\ w\in \Sigma^*\}$, it is useless and can be
removed. From now on we assume that
$G = \{\init\act w\ :\ w\in \Sigma^*\}$. \medskip

Of course, since $w\act{\mathscr M} = {\mathscr M}\act\overline{w}$
(where $\overline{w}$ stands for the mirror word of~$w$) it would be
enough to consider right actions only. But using both actions will
prove convenient in what follows.

We say that two Moore machines ${\mathscr M}$ and ${\mathscr M}'$,
with the same~$\Sigma$ and~$\Delta$, are equivalent if, for all
$w\in \Sigma^*$, ${\mathscr M}\act w = {\mathscr M}'\act w$.

Let ${\mathscr M}$ and ${\mathscr M}'$ be two Moore machines. If there
exists a bijection~$\xi$ from~$Q$ onto~$Q'$ such that, for all
$a\in Q$ and $j\in \Sigma$, $\delta'(\xi(a),j) = \xi(\delta(a),j)$
and ${\lambda}'\bigl(\xi(a)\bigr)= {\lambda}(a)$, we say that these
two machines are isomorphic.

\subsection{Duals}

\subsubsection*{Actions of words on functions}

If $f$ is a function from $Q$ to some set and $w\in \Sigma^*$ we
define two functions $ w\act f$ and $f\act w$:

$$\text{for all\quad } a\in Q,\quad (w\act f)(a) = f(a\act w)  \text{\quad and\quad}
(f\act w)(a) = f(w\act a)$$\medskip

\noindent Then\quad
$\bigl((vw)\act f\bigr)(a) = f\bigl( a\act(vw)\bigr) = (w\act f)(a\act
v) = \bigl(v\act(w\act f)\bigr)(a)$. So,
$$(vw)\act f = v\act(w\act f)$$ and similarly
$$f\act(vw) = (f\act v)\act w$$

\subsubsection*{Right dual} This is the Moore machine
${\mathscr M}_{\star} =
(Q_{\star},\Sigma,\Delta,\delta_{\star},{\lambda}_{\star},\init_{\star})$
defined as follows:\\
$Q_{\star} = \{w\act \lambda\ :\ w\in \Sigma^*\} \subset \Delta^Q$\\
$\delta_{\star}(f,j) = j\act f$\\
${\lambda}_{\star}\ :\ Q_{\star}\longrightarrow \Delta$ so defined:
\quad ${\lambda}_{\star}(f)= f(\init)$\\
$\init_{\star} = {\lambda}$.\\
\medskip

Then, for all $w\in \Sigma^*$,
${\lambda}_{\star}(w\act \init_{\star}) = (w\act \init_{\star})(\init) = (w\act
{\lambda})(\init) = {\lambda}(\init\act w)$, which means
  $$w\act {\mathscr M}_{\star} = {\mathscr M}\act w$$

  \subsubsection*{Left dual} We define the machine
  $\leftstar{\mathscr M} =
  (\leftstar{Q},\Sigma,\Delta,\leftstar{\delta},\leftstar{{\lambda}},
  \leftstar\init)$ in the same way:\\
  $\leftstar{Q} = \{ {\lambda}\act w\ :\ w\in \Sigma^*\} \subset \Delta^Q$\\
  $\leftstar{\delta}(f,j) = f\act j$\\
  $\leftstar{{\lambda}}\ :\ Q_{\star}\longrightarrow \Delta$
  so defined:\quad ${\lambda}_{\star}(f)= f(\init)$,\\
  $\leftstar{\init} = {\lambda}$. \medskip

  Then, for all $w\in \Sigma^*$,
  $\leftstar{{\lambda}}(\leftstar{\init}\act w) = (\leftstar{\init}\act w)(\init) = ({\lambda}\act
  w)(\init) = {\lambda}(w\act \init)$, which means
$$\leftstar{\mathscr M}\act w = w\act{\mathscr M}$$

\subsection{Bidual}

To avoid cumbersome notation, set
${\mathscr M}^\natural = \leftstar{({\mathscr M}_{\star})}$. Then both
$\mathscr M$ and ${\mathscr M}^\natural$ are right machines and are
equivalent in that sense that they give the same outputs.\medskip

There is a natural mapping~$\tau$ from $Q$ into $\Delta^{Q_{\star}}$:\\
let $a\in Q$ and $f\in Q_{\star}$, then $f\in \Delta^Q$, so
$f(a)\in \Delta$; set $\tau_a(f) = f(a)$. Then
$\tau\,:\,a\mapsto \tau_a$ is a mapping from~$Q$
to~$\Delta^{Q_{\star}}$.\medskip

We have $({\lambda}_{\star}\act w)(f) = {\lambda}_{\star}(w\act f) =
(w\act f)(\init) =
f(\init\act w) = \tau_{\init\cdot w}(f)$\\
for all $w\in \Sigma^*$ and $f\in Q_{\star}$. So
$${\lambda}_{\star}\act w = \tau_{\init\cdot w} \text{\quad for all\quad} w\in
\Sigma^*.$$
Since
$Q^\natural = \{\lambda_{\star}\act w\ :\ w\in \Sigma^*\} \subset
\Delta^{Q_{\star}}$, then $\tau$ maps $Q$ \textbf{onto} ${Q}^\natural$,
and $\card {Q}^\natural \le \card Q$.\medskip

One would like to know when $\card {Q}^\natural=\card Q$, that is
when~$\tau$ is one-to-one.  Let $a\in Q$ and $b\in Q$. Then
\begin{eqnarray*}
\tau_a=\tau_b &\Longleftrightarrow& \forall f\in Q_{\star}, f(a)=f(b)\\
              &\Longleftrightarrow& \forall w\in \Sigma^*, (w\act \lambda)(a)=(w\act \lambda)( b)\\
  &\Longleftrightarrow& \forall w\in \Sigma^*, {\lambda}(a\act w) =
                      {\lambda}(b\act w).
\end{eqnarray*}

Also, one can consider the mapping $\tau_{\star}$ from $Q_{\star}$ to
$\Delta^{{Q}^\natural}$ defined in the same way as $\tau$ is:
$$\tau_{\star,a}(f) = f(a) \text{\quad for \quad} a\in Q\_\star \text{ and } f\in Q^\natural.$$
Let $a_{\star}\in Q_{\star}$ and $b_{\star}\in Q_{\star}$. Then
\begin{eqnarray*}
  \tau_{\star,a_\star}=\tau_{\star,b_\star} &\Longleftrightarrow&  \forall f\in Q^\natural, f(a_\star) = f(b_\star)\\
                                            &\Longleftrightarrow& \forall w\in \Sigma^*, (\lambda_\star\act w)(a_\star) = (\lambda_\star\act w)(b_\star)\\
                                            &\Longleftrightarrow&  \forall w\in \Sigma^*, \lambda_\star(w\act a_\star) = \lambda_\star(w\act b_\star)\\ 
                                            &\Longleftrightarrow&  \forall w\in \Sigma^*, (w\act a_\star)(\init) = (w\act b_\star)(\init)\\
                                            &\Longleftrightarrow& \forall w\in \Sigma^*, a_{\star}(\init\act w) =b_{\star}(\init\act w)\\ &\Longleftrightarrow& \forall a\in Q, a_{\star}(a)=b_{\star}(a)\\
                                            &\Longleftrightarrow& a_\star = b_\star.
\end{eqnarray*}

This means that $\tau_{\star}$ is one-to-one and that
${\mathscr M}_{\star}$ is equal to its bidual.  In the same way
$\leftstar{\mathscr M}$ is equal to its bidual. Therefore
${\mathscr M}^\natural$ is equal to its bidual.

  \subsection{The algorithm}\label{algo}

  We illustrate  these constructions with the following example.
 \hspace*{-7em}

 \begin{center}
   \raisebox{-4em}{\includegraphics[width=70mm]{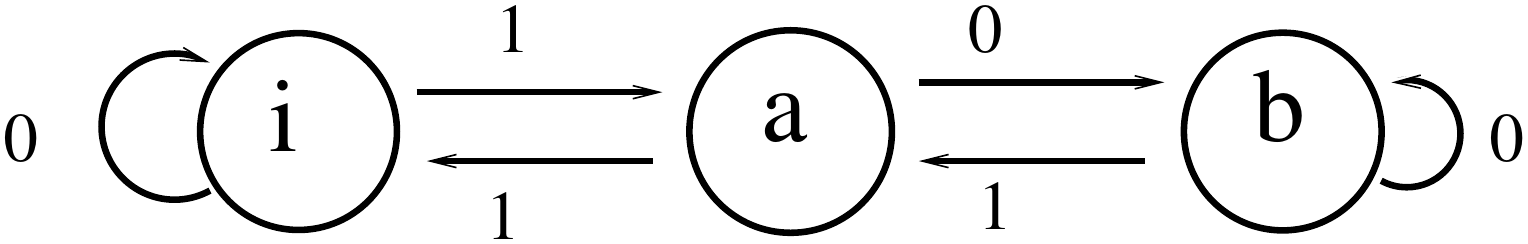}}\qquad
  $\begin{array}[t]{c||c|c}\multicolumn{3}{c}{ \text{Transition}}\\\hline
    &0&1\\\hline   \init&\init&a\\a&b&\init\\b&b&a\end{array}\quad
  \begin{array}[t]{c|c}\multicolumn{2}{c}{\lambda}\\\hline
    \init&0\\a&1\\b&0\end{array}$
\end{center}
\bigskip

Here $Q=\{\init,a,b\}$, $\Sigma = \{0,1\}$, and $\Delta =\{0,1\}$. The
transition and the output function are defined by the above
tables. Also, as it is customary, this machine is described by the
above diagram.

Now we describe an algorithm on this example. We consider a stack
whose members are elements of~$\Delta^Q$. We start with the stack
whose~$\lambda$ is the only element.We say that an element~$\alpha$ of
the stack is happy if $0\act\alpha$ and $1\act\alpha$ are both
elements of the stack. To add elements to the stack proceed as
follows: starting from the bottom take the first unhappy element,
say~$\alpha$, put $0\act\alpha$ on the top, set
$\delta_\star(\alpha,0)=0\act\alpha$, and do the same
with~$1\act\alpha$.Repeat until all the elements are happy. Then
$Q_\star$ is the set of elements of the stack and the transition
function is~$\delta_\star$.

The following table is the result of this process. A column whose
first element is $w\act \lambda$ contains two sub-columns; the left
one contains the states $\init\act w,\,a\act w$, and $b\act w$, the right one
contains $w\act \lambda$, i.e.,
$\lambda(\init\act w),\,\lambda(a\act w)$, and $\lambda(b\act w)$. The last
line is just a renaming of the states. This is illustrated by the
diagram and the tables below.\bigskip
  
 $\begin{array}{lr|lr|lr|lr|lr|lr|lr|lr|lr}
    \multicolumn{2}{c|}{\lambda}& \mult{0}{\lambda}& \mult{1}{\lambda}&\mult{00}{\lambda}& \mult{10}{\lambda}& \mult{01}{\lambda}& \mult{11}{\lambda}& \mult{001}{\lambda}&\multicolumn{2}{c}{101\act \lambda}\\
    \hline\hline
    \init&0&\init&0&a&1&\init&0&b&0&a&1&\init&0&a&1&a&1\\ \hline
    a&1&b&0&\init&0&b&0&\init&0&a&1&a&1&a&1&a&1\\ \hline
    b&0&b&0&a&1&b&0&b&0&a&1&\init&0&a&1&a&1\\ \hline\hline
    \multicolumn{2}{c}{t}& \multicolumn{2}{c}{u}& \multicolumn{2}{c}{v}& \multicolumn{2}{c}{u}& \multicolumn{2}{c}{u}& \multicolumn{2}{c}{w}& \multicolumn{2}{c}{t}& \multicolumn{2}{c}{w}& \multicolumn{2}{c}{w}
  \end{array}$

  \includegraphics[width=50mm]{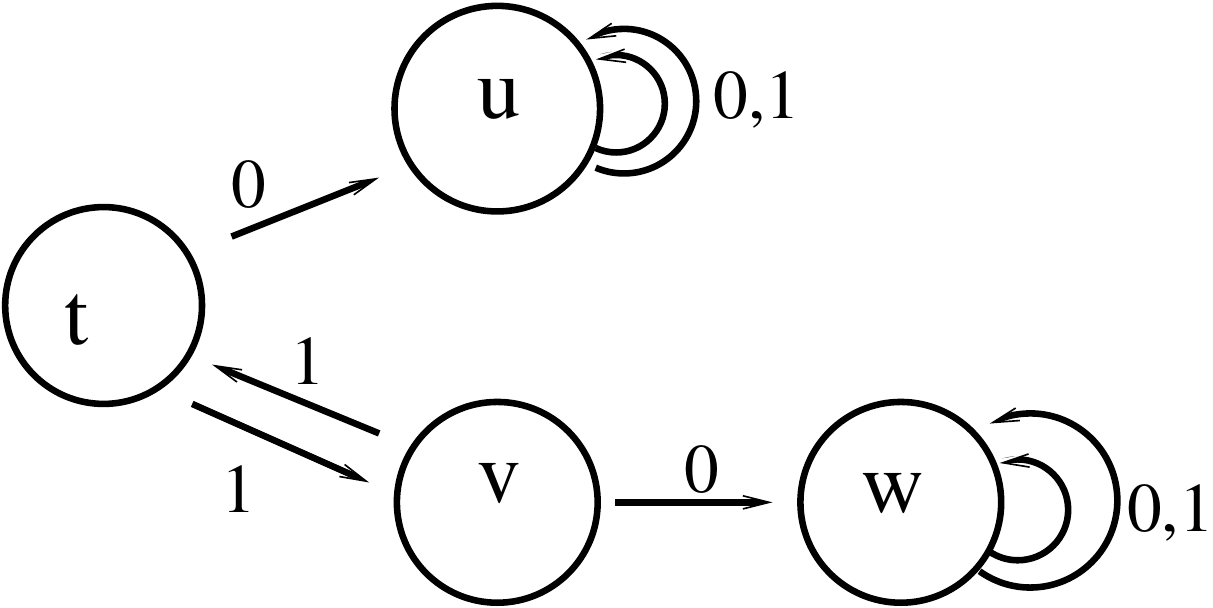}\qquad 
$  \begin{array}[b]{c||c|c}\multicolumn{3}{c}{ \text{Transition}}\\\hline
    &0&1\\\hline   t&u&v\\u&u&u\\v&w&t\\w&w&w\end{array}\quad
  \begin{array}[b]{c|c}\multicolumn{2}{c}{\lambda_{\star}}\\\hline
    t&0\\u&0\\v&1\\w&1\end{array}$
  \bigskip
  
  The following is the computation of the bidual (the states $u$ and
  $w$ have been omitted in the next table because they would not have
  brought any information). \medskip

$  \begin{array}{|lr|lr|lr|lr|lr|}
    \multicolumn{2}{c}{\lambda_{\star}}& \multicolumn{2}{c}{\lambda_{\star}\act 0}& \multicolumn{2}{c}{\lambda_{\star}\act 1}&\multicolumn{2}{c}{\lambda_{\star}\act 10}& \multicolumn{2}{c}{\lambda_{\star}\act 11}\\
    \hline\hline
    t&0&u&0&v&1&u&0&t&0\\ \hline
    v&1&w&1&t&0&w&1&v&1\\ \hline\hline
    \multicolumn{2}{c}{x}& \multicolumn{2}{c}{x}& \multicolumn{2}{c}{y}& \multicolumn{2}{c}{x}& \multicolumn{2}{c}{x}
  \end{array}$

  \includegraphics[width=40mm]{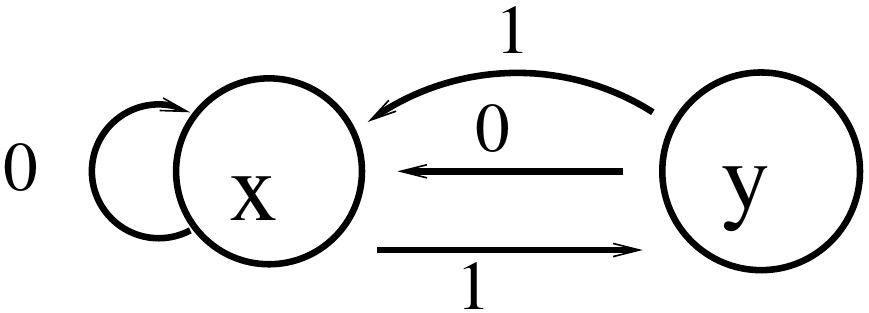}\hspace{5em}
$  \begin{array}[b]{c||c|c}\multicolumn{3}{c}{ \text{Transition}}\\\hline
    &0&1\\\hline   x&x&y\\y&x&x\end{array}\quad
  \begin{array}[b]{c|c}\multicolumn{2}{c}{{\lambda}^\natural}\\\hline
    x&0\\y&1\end{array}$

As expected, the bidual has less states as the original machine.

\section{Minimal Moore machine}

Here we give another proof of the existence and uniqueness of the
minimal machine.

\subsection{Product of machines}

Let
${\mathscr M}_1= (Q_1,\Sigma, \Delta_1,\delta_1,\lambda_1,\init_1)$
and
${\mathscr M}_2= (Q_2,\Sigma, \Delta_2,\delta_2,\lambda_2,\init_2)$
be two Moore machines, $\Delta$ a finite set, and $\gamma$ a mapping
from $\Delta_1\times \Delta_2$ to~$\Delta$. We define a  new machine
$${\mathscr M}_1\otimes_{\gamma}{\mathscr M}_2 =
\bigl(Q,\Sigma,\Delta,\delta,\gamma\circ(\lambda_1\times
\lambda_2),(\init_1,\init_2)\bigr)$$
by setting
$$Q = \{
(\init_1\act w,\init_2\act w)\ : w\in \Sigma^* \} \text{ and }
\delta\bigl( (a_1,a_2),j\bigr) = \bigl(
\delta_1(a_1,j),\delta_2(a_2,j)\bigr).$$

\subsection{Minimal Moore machine}

Let ${\mathscr M}_1$ and ${\mathscr M}_2$ be two equivalent right
machines such that each one is isomorphic to its bidual.\medskip

Let $p_1$ and $p_2$ be the projections of $Q_1\times Q_2$ onto $Q_1$
and $Q_2$. Consider the machines
${\mathscr B}_k = {\mathscr M_1}\otimes_{p_k} {\mathscr M}_2$, for
$k=1,2$.\medskip

For $w\in \sigma^*$, $k\in \{1,2\}$, and $(a_1,a_2)\in Q_1\times Q_2$ we have
\begin{eqnarray*}
  \Bigl(w\act\bigl(p_k\circ(\lambda_1\times\lambda_2)\bigr)\Bigg)(a_1,a_2) &=&  p_k\circ(\lambda_1\times \lambda_2)(a_1\act w, a_2\act w)\\
                                                                            &=& \lambda_k(a_k\act w) = (w\act\lambda_k)(a_k).
\end{eqnarray*}
This means
$w\act \bigl(p_k\circ(\lambda_1\times\lambda_2)\bigr) =
(w\act\lambda_y)\circ p_k$, which implies
${\mathscr B}_{k\star} = {\mathscr M}_{k\star}$.  Therefore
${\mathscr B}_k^\natural = {\mathscr M}_k^\natural = {\mathscr
  M}_k$\medskip

Because ${\mathscr M}_1$ and ${\mathscr M}_2$ are equivalent we have
$\lambda_1(a_1) = \lambda_2(a_2)$ for all $(a_1,a_2)\in Q$ (i.e.,
$a_1=\init_1\act w$ and $a_2=\init_2\act w$ for some $w\in
\Sigma^*$). Therefore ${\mathscr B}_1 = {\mathscr B}_2$.\medskip

As a conclusion, for any Moore machine, there is a unique (up to
isomorphism) simplest equivalent machine which is the bidual. One
could notice that this algorithm provides a normal form for a Moore
machine.
  
\section{Substitutions}


Consider a 5-tuple
${\mathscr S} = (Q,\Delta,\sigma,\lambda,\init)$ where
\begin{itemize}
\item[--] $Q$ and $\Delta$ are finite sets,
\item[--] $\sigma$ is an endomorphism of the free monoid~$Q^*$
  generated by~$Q$,
\item[--] $\lambda$ is a mapping from $Q$ to~$\Delta$,
\item[--] $\init\in Q$.
\end{itemize}

The endomorphism~$\sigma$, also called substitution, is determined by
the images of the generators, i.e., by the
words~$\{\sigma(a)\}_{a\in Q}$.  We adopt the following notation: if
$w$ is a word on some alphabet, $|w|$ stands for its length (i.e., the
number of letters it is made of), and $w_j$ stands for its $(j+1)$-th
letter from the left (the letters are numbered 0,\,1,\,2,\,\dots).
Let $q=\displaystyle \max_{a\in Q}|\sigma(a)|$.

When $\delta(\init,0)=\init$, the words
$\bigl(\sigma^n(\init)\bigr)_{n\ge 0}$ have a limit
$\sigma^\infty(\init)$ in $Q^{\mathbb N}$. Such sequences as
$\lambda\bigl( \sigma^\infty(\init)\bigr)$, called projections of
fixed points of substitutions, have been studied in many contexts,
namely for their algebraic, dynamic, and combinatoric
properties~\cite{allouche,christol,berthe}.

\subsection*{Constant length substitutions}
This is the cas when $|\sigma(a)|=q$ for all~$a\in Q$. We then set
$\Sigma= \{0,1,\dots,q-1\}$ and $\delta(a,j) = \sigma(a)_j$ for
$a\in Q$ and $0\le j< q$. So we have a Moore machine:
$${\mathscr M} =(Q,\Sigma,\Delta,\delta,\lambda,\init).$$
In this context if $k\ge 0$ and $n = n_0+n_1q+\dots n_{k-1}q^{k-1}$
(the base~$q$ expansion of~$n$) it is easy to see that the letter of
index~$n$ in the word $\sigma^k(a)$ is $n_0n_1\dots n_{k-1}\act a$,
i.e.,
\begin{equation} \label{enum1} \bigl(\sigma^k(a)\bigr)_n = n_0n_1\dots
  n_{k-1}\act a.
\end{equation}

This is the connection between Moore machines and substitutions. The
bidual of the left Moore machine~$\mathscr M$ defines a substitution
on a smaller alphabet equivalent to the original one.

\subsection*{Non constant length substitutions}

Let~$\omega$ be a symbol not in~$Q$. Set $Q'=Q\cup\{\omega\}$ and
$\Sigma=\{0,1,\dots,q-1\}$. For $0\le j< q$ set
$\delta(\omega,j)=\omega$, and for $a\in Q$ set
$\delta(a,j) = \bigl(\sigma(a)\omega^\nu\bigr)_j$, where
$\nu = q-|\sigma(a)|$.

We just have defined a finite
automaton~$\mathscr A = (Q',\Sigma,\delta,\init,Q)$, where~$\init$ is
the initial state and $Q$ is the set of final states.  Let
$\mathscr L$ be the language recognized by~$\mathscr A$ considered as
a left automaton :
${\mathscr L} = \{w\in \Sigma^*\ :\ w\act\init \ne \omega\}$. If we
define $\lambda(\omega)$ arbitrarily, we also have a Moore machine,
${\mathscr M} = (Q',\Sigma,\Delta,\delta,\lambda,\init)$. Let
${\mathscr L}' = {\mathscr L}\setminus\{\epsilon\}$. Observe that, if
$w\in {\mathscr L}'$, the path leading to $w\act \init$ does not
include any~$\omega$.

If $w\in {\mathscr L}'$, set
$\varphi(w) = \displaystyle \sum_{0\le j< |w|} w_j q^j$. Then the
order on $\mathbb N$ induces a well-order on ${\mathscr L}'$:
$v\prec w \Leftrightarrow \varphi(v)< \varphi(w)$. So, we can define
$\psi (n)$ to be the $n$-th element of ${\mathscr L}'$ (starting
from~$0$). Then the counterpart of~\eqref{enum1} is
\begin{equation*}
\bigl(\sigma^k(\init)\bigr)_j = \psi(j)\act\init
\end{equation*}
for $k\ge 0$ and $0\le j< |\sigma^k(\init)|$.

As previously, the bidual of the left Moore machine~$\mathscr M$
defines a substitution on a smaller alphabet equivalent to the
original one.

Nevertheless, observe that instead of appending some~$\omega$'s one
could have padded the $\sigma(a)$'s up to size~$q$ by
inserting~$\omega$'s in various places. Then, of course, the
language~$\mathscr L$ would be different, but the above analysis still
holds.

\end{document}